\documentclass[reqno]{amsart}
\usepackage{a4wide}

\theoremstyle{plain}
\newtheorem{theorem}{Theorem}
\newtheorem{corollary}{Corollary}

\newtheorem{lemma}{Lemma}

\newtheorem*{theorem*}{Theorem}
\newtheorem*{corollary*}{Corollary}
\newtheorem*{proposition*}{Proposition}
\newtheorem*{lemma*}{Lemma}

\theoremstyle{definition}

\newtheorem*{ack}{Acknowledgements}

\theoremstyle{remark}

\newcommand\re{\mathrm{Re}\,}
\newcommand\im{\mathrm{Im}\,}
\newcommand\bigsqrt[1]{\sqrt{\smash{#1} \vphantom{|^|}}}

\begin{document}
\title[Schr{\"o}dinger operators with $L^2$-sparse potentials]
{Value distribution and spectral theory of Schr{\"o}dinger operators with $L^2$-sparse potentials}
\author{S.V.\ Breimesser}
\author{James D.E.\ Grant}
\author{D.B.\ Pearson}
\address{Department of Mathematics\\ University of Hull\\ Hull HU6 7RX\\ U.K.}
\email{s.v.breimesser@hull.ac.uk}
\email{j.d.grant@hull.ac.uk}
\email{d.b.pearson@hull.ac.uk}
\subjclass[2000]{47E05, 34L05, 82Q10}
\keywords{Herglotz functions, Schr{\"o}dinger operator, value distribution, spectral theory, sparse potentials}
\begin{abstract}
We apply the methods of value distribution theory to the 
spectral asymptotics of Schr{\"o}dinger operators with 
$L^2$-sparse potentials.
\end{abstract}

\maketitle
\thispagestyle{empty}

\section{Introduction}

A real valued, locally integrable function $V$, defined on the half-line $0 \le x < \infty$, 
is said to be a $L^2$-sparse potential if, given arbitrary $\delta, N > 0$, 
there exists a subinterval $\left( a, b \right)$ of $\left[ 0, \infty \right)$ 
such that $b - a = N$ and $\int_a^b \left( V(x) \right)^2 dx < \delta$. 
In other words, if $V$ is $L^2$-sparse then one can find arbitrarily long 
intervals on which the $L^2$ norm of $V$ is arbitrarily small. 
Given an $L^2$-sparse potential, we can define a Schr{\"o}dinger operator 
$T = - \frac{d^2}{dx^2} + V(x)$ acting in $L^2 (0, \infty)$ 
and subject to Dirichlet boundary condition at $x=0$. 
By considering an appropriate sequence of approximate eigenfunctions 
(see for example \cite{G}, Theorem~22) one may verify that the Weyl spectrum 
of $T$ contains the whole of $\mathbb{R}^+$. It follows that we have the 
limit point case at infinity, so that $T$ can be uniquely defined as a 
self-adjoint operator, subject to the single boundary condition at $x=0$.

Any $L^2$-sparse potential is a sum $V_1 + V_2$, 
where $V_1$ is a sparse potential and $V_2 \in L^2 \left( 0, \infty \right)$; 
here a potential $V$ is said to be sparse if arbitrarily long intervals 
exist on which $V$ is identically zero. 
There is a considerable literature on sparse potentials and their perturbations, 
in particular establishing conditions for the existence of 
absolutely continuous and singular continuous spectra. 
For recent results in this field, see \cite{KLS,R,SS} and references therein.

Spectral theory for the Schr{\"o}dinger operator $T$ can be closely linked 
to the theory of value distribution for real-valued functions, 
and in particular value distribution for functions which are 
defined as boundary values of Herglotz functions. 
(A Herglotz or Nevanlinna function is a function of a complex 
variable, analytic in the upper half-plane with positive 
imaginary part.) 

For a measurable function $F_+: \mathbb{R} \rightarrow \mathbb{R}$, 
the value distribution may be described by means of a map   
$\mathcal{M}: \left( A, S \right) \mapsto 
\mathcal{M} \left( A, S \right) \in \mathbb{R} \cup \left\{ \infty \right\}$, 
called the value distribution function of $F_+$, 
and defined for Borel subsets $A, S$ of $\mathbb{R}$ by
\begin{equation}
\mathcal{M} \left( A, S \right) = \left| A \cap F_+^{-1}(S) \right|.
\label{valuedistn}
\end{equation}
Here $\left| \, \cdot \vphantom{|^|} \, \right|$ stands for Lebesgue measure. 
Thus $\mathcal{M} \left( A, S \right)$ is the Lebesgue measure of the 
set of $\lambda \in A$ for which $F_+(\lambda) \in S$. 
In the particular case that $F_+$ is the almost everywhere boundary 
value of a Herglotz function, i.e. 
\[
F_+(\lambda) = \lim_{d~\rightarrow~0^+} F(\lambda + i d), \qquad
\mbox{(almost all $\lambda \in \mathbb{R}$)},
\]
we can write (see \cite{BP1})
\begin{equation}
\mathcal{M} \left( A, S \right) = \lim_{d~\rightarrow~0^+} 
\frac{1}{\pi} \int_A \theta \left( F(\lambda + i d), S \right) d\lambda,
\label{vdint}
\end{equation}
where $\theta(z, S)$ denotes the angle subtended at a point $z$ by 
the Borel subset $S$ of the real line. (For $\lambda \in \mathbb{R}$, 
we define $\theta(\lambda, S)$ to be $\pi \chi_S(\lambda)$, where 
$\chi_S$ is the characteristic function of the set $S$.) 
In fact, given $A$ and $S$ with $\left| A \right| < \infty$, 
the limit in~\eqref{vdint} will exist for \emph{any} Herglotz 
function $F$ (whether or not $F$ has real boundary values a.e.) 
and may be used to define the value distribution function $\mathcal{M}$ 
associated with an arbitrary Herglotz function. In general $\mathcal{M}$ 
may not describe the value distribution of any single real-valued function 
$F_+(\lambda)$, but there will always be sequences $\{ F^{(n)} \}$ of 
real-valued functions for which $\mathcal{M}$ describes the limiting value 
distribution.

Value distribution for boundary values of Herglotz functions is also closely 
connected with the geometric properties of the upper half-plane, regarded 
as a hyperbolic space \cite{BP1,BP2}. Given two points $z_1, z_2 \in \mathbb{C}^+$, 
we define a measure of separation 
\begin{equation}
\gamma \left( z_1, z_2 \right) = 
\frac{\left| z_1 - z_2 \right|}{\bigsqrt{\im z_1} \bigsqrt{\im z_2}},
\label{gammadef}
\end{equation}
which is related to the hyperbolic distance $D(z_1, z_2)$ by the equation
\[
\gamma \left( z_1, z_2 \right) = 2 \sinh \left( \frac{1}{2} D(z_1, z_2) \right).
\]
The relevance of hyperbolic distance to estimates of value distribution 
comes from the fact that if $F_1$ and $F_2$ are two Herglotz functions 
satisfying the estimate
\[
\gamma \left( F_1(z), F_2(z) \right) < \epsilon,
\]
for all $z$ such that $\im z = d$ and $\re z \in A$, 
then the value distribution $\mathcal{M}_2(A, S)$ 
associated with $F_2$ is a good approximation to the value distribution 
$\mathcal{M}_1(A, S)$ associated with $F_1$, in the sense that
\begin{equation}
\left| \mathcal{M}_1(A, S) - \mathcal{M}_2(A, S) \vphantom{|^|} \right| 
\le \epsilon \left| A \right| + 2 E_A(d).
\label{vddiffest}
\end{equation}
Here $E_A(d)$ is an error estimate which is an increasing function of $d$, 
and which converges to zero in the limit $d \rightarrow 0$, 
for fixed Borel set $A$. For details of this and related results, see \cite{BP1,BP2}. 
Estimates such as~\eqref{vddiffest} imply that if 
$\{ F^{(n)} \}_{n=1, 2, 3, \dots}$ is a sequence of Herglotz functions 
converging uniformly to $F(z)$ for $z$ lying in any fixed compact 
subset of $\mathbb{C}^+$, then the value distribution associated with 
$F^{(n)}$ will converge in the limit $n \rightarrow \infty$ to the value 
distribution associated with $F$.

The main purpose of this paper is to apply the above analysis to the 
spectral asymptotics of the Schr{\"o}dinger operator with $L^2$-sparse potential, 
as described by solutions $f(x, \lambda)$ of the Schr{\"o}dinger equation at 
real spectral parameter $\lambda$. Herglotz functions of particular 
interest in this context are the Weyl $m$-function $m(z)$ for the 
operator $- \frac{d^2}{dx^2} + V$ in $L^2 \left( 0, \infty \right)$ 
with Dirichlet boundary condition at $x=0$, and the Weyl $m$-function 
$m^N(z)$ for $- \frac{d^2}{dx^2} + V$ regarded as operating in 
$L^2 \left( N, \infty \right)$ for some fixed $N>0$, 
with Dirichlet boundary condition at $x=N$. Estimates of both of these 
$m$-functions may be carried out, for $z$ in some fixed compact subset of 
$\mathbb{C}^+$, in terms of the logarithmic derivative 
$f^{\prime}(x, z) / f(x, z)$, for asymptotically large $x$, 
of particular solutions $f(\cdot, z)$ of the Schr{\"o}dinger equation at 
\emph{complex} spectral parameter $z$. The main general results 
of the paper are presented in Theorems~\ref{thm:1} and~\ref{thm:2}. 
Theorem~\ref{thm:1} provides an estimate of the large $x$ asymptotics 
of $f^{\prime}/f$, based on an $L^2$ bound for the potential across 
a finite interval. Theorem~\ref{thm:2} is an analysis of asymptotic 
value distribution in the case of $L^2$-sparse potential, linking 
this to the asymptotics of $m^N$. 

Finally, we indicate some consequences of the 
analysis for spectral theory of $L^2$-sparse potentials, 
implying in particular the absence of absolutely continuous spectrum at negative $\lambda$.

\section{Asymptotics of $v^{\prime} / v$}
We consider the differential expression $\tau = - \frac{d^2}{dx^2} + V(x)$ on the 
half-line $0 \le x < \infty$, where the potential function $V$ is assumed to be 
real valued and integrable over any finite subinterval of $\left[ 0, \infty \right)$. 
Assume limit-point case at infinity, implying that a self-adjoint operator 
$T = - \frac{d^2}{dx^2} + V(x)$ can be defined, acting in $L^2 (0, \infty)$ 
and subject to a Dirichlet boundary condition at $x=0$.

The Weyl $m$-function $m(z; V)$ may be defined in terms of solutions $f( \cdot, z)$ 
of the Schr{\"o}dinger equation at complex spectral parameter $z$, namely
\begin{equation}
- \frac{d^2}{dx^2} f(x, z) + V(x) f(x, z) = z f(x, z) \qquad 
(\im z > 0, 0 \le x < \infty).
\label{Schrodinger}
\end{equation}
First define two solutions $u(x, z), v(x, z)$ of~\eqref{Schrodinger}, 
subject respectively to initial conditions
\begin{subequations}
\begin{align}
u(0, z) = 1, & &u^{\prime}(0, z) = 0,
\label{uBC}
\\
v(0, z) = 0, & &v^{\prime}(0, z) = 1,
\label{vBC}
\end{align}
\label{BCs}
\end{subequations}
where prime denotes differentiation with respect to $x$. 
(Solutions of~\eqref{Schrodinger} and~\eqref{BCs} with $z$ replaced by a 
real spectral parameter $\lambda$ will be denoted by $u(x, \lambda), v(x, \lambda)$ 
respectively and, for fixed $x$, are the boundary values of 
$u(x, z)$, $v(x, z)$ as $z$ approaches the real axis.)

Then (in the limit-point case at infinity) we define $m(z; V)$ uniquely by the 
condition that
\[
u(\cdot, z) + m(z; V) \, v(\cdot, z) \in L^2 (0, \infty).
\] 
An alternative characterisation of the Weyl function is that if $f(\cdot, z)$ 
is any (non-trivial) $L^2 (0, \infty)$ solution of~\eqref{Schrodinger}, then 
\begin{equation}
m(z; V) = \frac{f^{\prime}(0, z)}{f(0, z)}.
\label{mdef}
\end{equation}
It follows from the limit point/limit circle theory \cite{CL} 
that $m(z; V)$ is an analytic function of $z$ for $\im z > 0$. 
In addition $\im m(z; V) > 0$ for $\im z > 0$, 
so that $m(z; V)$ is a Herglotz function 
(analytic in the upper half-plane with positive imaginary part).

Given any $N>0$, we can also define the Dirichlet $m$-function $m^N(z; V)$ for the 
truncated problem on the interval $N \le x < \infty$, and an analysis of the large 
$N$ asymptotics of $m^N$ will play an important role in this paper. Here we are 
strongly motivated by the recent results of Deift and Killip \cite{DK} for $L^2$ potentials. 
 
Since, according to equation~\eqref{mdef}, the $m$-function is 
dependent on the logarithmic derivative of a solution $f(\cdot, z)$ of 
equation~\eqref{Schrodinger}, a first step in our analysis will be to carry out 
a comparison between logarithmic derivatives of solutions of equation~\eqref{Schrodinger} 
as the potential is varied. We begin with the logarithmic derivative of the solution 
$v(\cdot, z)$ subject to the initial conditions~\eqref{vBC}. Here it is $- v^{\prime} / v$ 
rather than $v^{\prime}/v$ that is a Herglotz function for $x > 0$. 
The following elementary estimate provides a bound for the $\gamma$-separation of 
the logarithmic derivative as the potential is varied.
\begin{lemma}
\label{lem:gammaestimate}
Let $v(x, z), {\tilde v}(x, z)$ be solutions of equation~\eqref{Schrodinger} with 
potentials $V(x)$, ${\tilde V}(x)$ respectively, and subject to initial conditions 
\[
v(0, z) = {\tilde v}(0, z) = 0, \qquad
v^{\prime}(0, z) = {\tilde v}^{\prime}(0, z) = 1.
\]
Then, for any $x>0$, 
\begin{equation}
\gamma \left( - \frac{v^{\prime}(x, z)}{v(x, z)}, 
- \frac{{\tilde v}^{\prime}(x, z)}{{\tilde v}(x, z)} \right) \le 
\left( \im z \right)^{-1} 
\frac{\left( \int_0^x \left( V(t) - {\tilde V}(t) \right)^2 
\left| {\tilde v}(t, z) \right|^2 \, dt \right)^{1/2}}
{\left( \int_0^x \left| {\tilde v}(t, z) \right|^2 \, dt \right)^{1/2}}
\label{gammaestimate}
\end{equation}
\end{lemma}
\begin{proof}
Abbreviating the notation for simplicity, we have
\begin{equation}
\gamma \left( - \frac{v^{\prime}}{v}, 
- \frac{{\tilde v}^{\prime}}{{\tilde v}} \right) = 
\frac{\left| \frac{v^{\prime}}{v} - \frac{{\tilde v}^{\prime}}{{\tilde v}} \right|}
{\sqrt{ \im \left( - \frac{v^{\prime}}{v} \right) \, 
\im \left( - \frac{{\tilde v}^{\prime}}{{\tilde v}} \right)}},
\label{gam1}
\end{equation}
where
\begin{equation}
\left| \frac{v^{\prime}}{v} - \frac{{\tilde v}^{\prime}}{{\tilde v}} \right| = 
\frac{\left| {\tilde v} v^{\prime} - v {\tilde v}^{\prime} \right|}
{\left| v {\tilde v} \right|}.
\label{gam2}
\end{equation}
Using the Schr{\"o}dinger equation $-v^{\prime\prime} + V v = z v$, and similarly for ${\tilde v}$, 
we have
\[
\frac{d}{dx} \left( {\tilde v} v^{\prime} - v {\tilde v}^{\prime} \right) = 
{\tilde v} v^{\prime\prime} - v {\tilde v}^{\prime\prime} = 
\left( V - {\tilde V} \right) v {\tilde v},
\] 
which, with the initial conditions, gives 
\[
{\tilde v} v^{\prime} - v {\tilde v}^{\prime} = 
\int_0^x \left( V(t) - {\tilde V}(t) \right) v(t) {\tilde v}(t) \, dt.
\]
We also have
\[
\im \left( - \frac{v^{\prime}}{v} \right) = 
\frac{1}{2i} \left( \frac{{\overline v}^\prime}{{\overline v}} - 
\frac{v^{\prime}}{v} \right) = 
\frac{1}{2i |v|^2} \left( v {\overline v}^{\prime} - {\overline v} v^{\prime} \right),
\]
which again on considering 
$\frac{d}{dx} \left( v {\overline v}^{\prime} - {\overline v} v^{\prime} \right)$ gives
\[
\im \left( - \frac{v^{\prime}}{v} \right) = 
\frac{\im z}{|v|^2} \int_0^x \left| v(t) \right|^2 \, dt,
\]
with a similar equation for ${\tilde v}$. 
Using~\eqref{gam1} and \eqref{gam2}, and substituting for 
$\im \left( - v^{\prime} / v \right)$, 
$\im \left( - {\tilde v}^{\prime} / {\tilde v} \right)$ 
and $\left( {\tilde v} v^{\prime} - v {\tilde v}^{\prime} \right)$ 
results in the bound
\[
\gamma \left( - \frac{v^{\prime}}{v}, 
- \frac{{\tilde v}^{\prime}}{{\tilde v}} \right) = 
\frac{\left| \int_0^x \left( V(t) - {\tilde V}(t) \right) v(t) {\tilde v}(t) \, dt \right|}
{\left( \im z \right) 
\left( \int_0^x \left| v(t) \right|^2 \, dt 
\int_0^x \left| {\tilde v}(t) \right|^2 \, dt \right)^{1/2}},
\]
from which~\eqref{gammaestimate} follows on applying Schwarz's inequality 
to the integral in the numerator.
\end{proof}

If both potentials $V, {\tilde V}$ are bounded, 
we can use the result of Lemma~\ref{lem:gammaestimate} 
to derive simple bounds for the separation $\gamma$ 
between the two logarithmic derivatives. 
For example we have, from~\eqref{gammaestimate}, for any $L > 0$,
\[
\left. \gamma \left( - \frac{v^{\prime}}{v}, 
- \frac{{\tilde v}^{\prime}}{{\tilde v}} \right) \right|_{x=L} \le 
\frac{1}{\im z} \sup_{t \in \left[ 0, L \right]} 
\left| V(t) - {\tilde V}(t) \right|.
\]
In particular, we see that any \emph{uniformly} convergent sequence 
$V_n$ of potentials will result in a corresponding sequence 
$- v_n^{\prime} / v_n$ which will converge uniformly in $\gamma$-separation 
(and hence also uniformly in the hyperbolic metric). 

We turn now to the case of a potential subject to an $L^2$-type condition, 
for which we take in the first instance the comparison potential to be 
${\tilde V}(x) = 0$. Let $v(x, z)$ be defined as before to be the solution 
of equation~\eqref{Schrodinger} with potential $V(x)$ and subject to 
$v(0, z) = 0, v^{\prime}(0, z) = 1$, and let $v_0(x, z)$ satisfy the equation
\[
- \frac{d^2 v_0(x, z)}{dx^2} = z v_0(x, z)
\]
with the same initial conditions. Again we take $\im z > 0$, and write 
$\sqrt{z} = a + i b$ with $a, b$ real and $a, b > 0$. 
An explicit expression for $v_0$ is then 
\[
v_0(x, z) = \left( 2 i \sqrt{z} \right)^{-1} 
\left( \mathrm{e}^{i x \sqrt{z}} - \mathrm{e}^{- i x \sqrt{z}} \right) 
= \left( 2 \left( b - i a \right) \right)^{-1} 
\left( \mathrm{e}^{-iax} \mathrm{e}^{bx} - \mathrm{e}^{iax} \mathrm{e}^{-bx} \right),
\]
so that
\[
\left| v_0(x, z) \right|^2 = \left( 2 \left( a^2 + b^2 \right) \right)^{-1} 
\left( \cosh 2 b x - \cos 2 a x \right),
\]
and, from~\eqref{gammaestimate}, we have
\begin{equation}
\left. \gamma \left( - \frac{v^{\prime}}{v}, 
- \frac{v_0^{\prime}}{v_0} \right) \right|_{x=L} 
\le 
\frac{
\left( \int_0^L V(t)^2 \left( \cosh 2 b t - \cos 2 a t \right) \, dt \right)^{1/2}}
{\im z \left( \int_0^L \left( \cosh 2 b t - \cos 2 a t \right) \, dt \right)^{1/2}}.
\label{v=0estimate}
\end{equation}
Here the integral in the numerator may be written
\begin{align}
- \int_0^L \left\{ \left( \cosh 2 b t - \cos 2 a t \right) 
\frac{d}{dt} \int_t^L V(s)^2 \, ds \right\} \, dt 
&= \int_0^L \left\{ \left( 2 b \sinh 2 b t + 2 a \sin 2 a t \right) 
\int_t^L V(s)^2 \, ds \right\} \, dt
\nonumber
\\
&\le \int_0^L \left\{ \left( 2 b \sinh 2 b t + 2  a  \right) 
\int_0^L V(s)^2 \, ds \right\} \, dt
\nonumber
\\
&= \left( 2 a L + \cosh 2 b L - 1 \right) \int_0^L V(s)^2 \, ds
\label{v=0ineq}
\end{align}

To complete the estimate of~\eqref{v=0estimate}, we need a lower bound for the 
denominator integral, which comes to 
\[
\frac{\sinh 2 b L}{2b} - \frac{\sin 2 a L}{2a}.
\]
We shall make the assumption $L \ge 1 / \sqrt{|z|}$. Such a condition, 
with $\sqrt{z} = a + i b$, implies that either 
$L \ge 1 / \left( \sqrt{2} a \right)$ or $L \ge 1 / \left( \sqrt{2} b \right)$. 
(If $L <1 / \left( \sqrt{2} a \right)$ and $L < 1 / \left( \sqrt{2} b \right)$ 
then $|z| = a^2 + b^2 < \frac{1}{2L^2} + \frac{1}{2L^2} = \frac{1}{L^2}$, 
which contradicts the assumption.) 

We consider the two possibilities in turn:

\vskip .2cm
\noindent{\fbox{Case 1 : $L \ge \frac{1}{\sqrt{2} a}$}}
\vskip .2cm
{}From the bound $\sinh x / x > 1$ for $x > 0$, we have 
\[
\frac{\sinh 2 b L}{2b} > L, 
\]
whereas 
\[
\left| \frac{\sin 2 a L}{2a} \right| \le \frac{1}{2a} \le \frac{L}{\sqrt{2}},
\] 
so that 
\[
\left| \frac{\sin 2 a L}{2a} \right| < \frac{1}{\sqrt{2}} \frac{\sinh 2 b L}{2b},
\] 
and it follows that
\begin{equation}
\frac{\sinh 2 b L}{2b} - \frac{\sin 2 a L}{2a} > 
\left( 1 - \frac{1}{\sqrt{2}} \right) \frac{\sinh 2 b L}{2b}.
\label{sinhineq1}
\end{equation}

\vskip .2cm
\noindent{\fbox{Case 2 : $L \ge \frac{1}{\sqrt{2} b}$}}
\vskip .2cm
Since the function $\sinh x / x$ is increasing for $x \ge 0$, 
we then have 
\[
\frac{\sinh 2 b L}{2b} \ge \frac{L \sinh \sqrt{2}}{\sqrt{2}},
\]
whereas
\[
\left| \frac{\sin 2 a L}{2 a} \right| < L.
\]
Hence in this case we find
\[
\left| \frac{\sin 2 a L}{2 a} \right| < \frac{\sqrt{2}}{\sinh \sqrt{2}} 
\left( \frac{\sinh 2 b L}{2b} \right),
\]
so that
\begin{equation}
\frac{\sinh 2 b L}{2b} - \frac{\sin 2 a L}{2a} > 
\left( 1 - \frac{\sqrt{2}}{\sinh \sqrt{2}} \right) \frac{\sinh 2 b L}{2b}.
\label{sinhineq2}
\end{equation}

Noting that $\sinh \sqrt{2} < 2$, we see that the bound~\eqref{sinhineq2} 
holds both in case~1 and in case~2.

\vskip .2cm

Using~\eqref{v=0ineq} and~\eqref{sinhineq2} as upper and lower bounds for the 
numerator and denominator respectively of~\eqref{v=0estimate}, we have, now, 
for $L \ge 1 / \sqrt{|z|}$, the estimate
\begin{align*}
\gamma \left( - \frac{v^{\prime}(L, z)}{v(L, z)}, 
- \frac{v_0^{\prime}(L, z)}{v_0(L, z)} \right) 
&\le 
\frac{1}{\im z} 
\left( 
\frac{\left( 2 a L + \cosh 2 b L - 1 \right) \int_0^L V(s)^2 \, ds}
{\left( 1 - \frac{\sqrt{2}}{\sinh \sqrt{2}} \right) \frac{\sinh 2 b L}{2b}}
\right)^{1/2}
\\
&= 
\frac{1}{\im z} \left( 1 - \frac{\sqrt{2}}{\sinh \sqrt{2}} \right)^{-1/2}
\left( 2 a \left( \frac{2bL}{\sinh 2 b L} \right) + 
2 b \left( \frac{\cosh 2 b L - 1}{\sinh 2 b L} \right) 
\right)^{1/2}
\\
& \hskip 4cm
\times \left( \int_0^L V(s)^2 \, ds \right)^{1/2}
\end{align*}
Noting that
\[
\frac{2bL}{\sinh 2bL} < 1
\]
and that
\[
\frac{\cosh 2 b L - 1}{\sinh 2 b L} = \tanh bL < 1,
\]
we can use the estimate 
$\left( a + b \right)^{1/2} \le 
\left( 2 \left( a^2 + b^2 \right) \right)^{1/4} = 
\left( 2 |z| \right)^{1/4}$ 
to obtain the following result.
\begin{lemma}
\label{lem:2}
Define $v(x, z)$ as in Lemma~\ref{lem:gammaestimate}, and let $v_0(x, z)$ be 
the corresponding solution of~\eqref{Schrodinger} with zero potential. 
Then, for any $L \ge 1 / \sqrt{|z|}$, we have the bound
\begin{equation}
\gamma \left( - \frac{v^{\prime}(L, z)}{v(L, z)}, 
- \frac{v_0^{\prime}(L, z)}{v_0(L, z)} \right) \le 
\frac{C |z|^{1/4}}{\im z} \left( \int_0^L V(s)^2 \, ds \right)^{1/2},
\label{z1/4estimate}
\end{equation}
where $C$ is a positive constant. 
(In fact we can take $C = \sqrt{2} \left( \frac{1}{\sqrt{2}} - \frac{1}{\sinh{\sqrt{2}}} \right)^{-1/2}$ 
in which case $C < 3.3$.)
\end{lemma}

Notice that Lemma~\ref{lem:2} provides a simple bound for the hyperbolic distance between 
$-v^{\prime}/v$ and $-v_0^{\prime}/v_0$ at $x=L$, in terms of the $L^2$ norm of the potential $V$ 
across the interval $\left[ 0, L \right]$. 

Since, as is easily verified, we have 
\[
\lim_{L~\rightarrow~\infty} - \frac{v_0^{\prime}(L, z)}{v_0(L, z)} = i \sqrt{z},
\] 
we can make a comparison, for large $L$, of $-v_0^{\prime}/v_0$ with its asymptotic 
limit, leading to the following result.

\begin{lemma}
\label{lem:3}
With $v_0(x, z)$ defined as in Lemma~\ref{lem:2}, for any $L \ge 1/\sqrt{\left| z \right|}$ 
we have the bound
\begin{equation}
\gamma \left( - \frac{v_0^{\prime}(L, z)}{v_0(L, z)}, i \sqrt{z} \right) \le 
C^{\prime} \frac{\left( 1 + \left( \frac{b}{a} \right)^2 \right)^{1/2}}
{\left( \mathrm{e}^{4bL} - 1 \right)^{1/2}},
\label{v0rootz}
\end{equation}
where $C^{\prime}$ is a positive constant. 
(In fact we can take $C^{\prime} = 2^{1/4} C$, where $C$ is the constant defined 
in Lemma~\ref{lem:2}, in which case $C^{\prime} < 3.9$.) 
\end{lemma}
\begin{proof}
Explicitly, we have
\[
- \frac{v_0^{\prime}}{v_0} = 
\frac{\left( i a - b \right) 
\left( \mathrm{e}^{-iax} \mathrm{e}^{bx} + \mathrm{e}^{iax} \mathrm{e}^{-bx} \right)}
{\mathrm{e}^{-iax} \mathrm{e}^{bx} - \mathrm{e}^{iax} \mathrm{e}^{-bx}},
\]
and multiplying numerator and denominator by the complex conjugate of the denominator gives
\[
\im \left( - \frac{v_0^{\prime}}{v_0} \right) = 
\frac{2 a \sinh 2 b x - 2 b \sin 2 a x}
{\left| \vphantom{|^|} \mathrm{e}^{-iax} \mathrm{e}^{bx} - \mathrm{e}^{iax} \mathrm{e}^{-bx} \right|^2},
\]
Moreover, 
\[
\left| - \frac{v_0^{\prime}}{v_0} - i \sqrt{z} \right| = 
\left| - \frac{v_0^{\prime}}{v_0} + b - i a \right| = 
\frac{2 \sqrt{a^2+b^2} \mathrm{e}^{-bx}}
{\left| \mathrm{e}^{-iax} \mathrm{e}^{bx} - \mathrm{e}^{iax} \mathrm{e}^{-bx} \right|}.
\]
Putting these results together we find, at $x=L$,
\[
\left. \gamma \left( - \frac{v_0^{\prime}}{v_0}, i \sqrt{z} \right) \right|_{x=L} = 
\left( \frac{2 \left( a^2 + b^2 \right)}{a} \right)^{1/2}
\frac{\mathrm{e}^{-bL}}{\left( a \sinh 2 b L - b \sin 2 a L \right)^{1/2}}.
\]
Substituting in the denominator the lower bound obtained previously in~\eqref{sinhineq2} 
and simplifying, we arrive at~\eqref{v0rootz}.
\end{proof}

In using~\eqref{v0rootz} to make precise estimates of the convergence to $i \sqrt{z}$ 
of $-v_0^{\prime}/v_0$, it is useful to note the inequalities:
\begin{itemize}
\item[(i)] if $\re z \ge 0$ then $b/a \le 1$;
\item[(ii)] $\left( 1 + \left( \frac{b}{a} \right)^2 \right) \le 
4 \left( 1 + \left( \frac{\re z}{\im z} \right)^2 \right)$; 
\item[(iii)] $b > \frac{\im z}{2 \sqrt{|z|}}$.
\end{itemize}
These inequalities imply, in particular, that $-v_0^{\prime}/v_0$ converges uniformly 
in hyperbolic norm to $i \sqrt{z}$, for $z$ in any fixed compact subset of the 
upper half-plane.

\section{Estimates of $u$ and $v$ for $L^1$-bounded potentials}
We consider solutions $u(x, z), v(x, z)$ of equation~\eqref{Schrodinger} on a fixed 
interval $0 \le x \le N$, subject to initial conditions~\eqref{BCs} at $x=0$. 
We compare these solutions with the corresponding solutions $u_0(x, z), v_0(x, z)$ 
with zero potential, and subject to the same initial conditions as for $u$ and $v$.

\begin{lemma}
\label{lem:4}
Let $K$ be a fixed compact subset of $\mathbb{C}^+$, 
and let $N>0$ be fixed. Then, given any $\epsilon > 0$, 
there exists $\delta_0 > 0$ ($\delta_0$ depending on 
$\epsilon, N$ and $K$) such that, for any potential function $V$ satisfying 
\[
\int_0^N \left| V(t) \right| \, dt < \delta_0,
\]
we have, for all $z \in K$ and for all $x \in \left[ 0, N \right]$,
\[
\left| u(x, z) - u_0(x, z) \right| < \epsilon, \qquad
\left| v(x, z) - v_0(x, z) \right| < \epsilon. 
\]
\end{lemma}
\begin{proof}
The proof is a standard perturbation argument using the Gronwall inequality. 

Let $M$ be the $2 \times 2$ transfer matrix given by
\[
M = M(x, z) = \begin{pmatrix} u & v \\ u^{\prime} & v^{\prime} \end{pmatrix},
\]
and let
\[
M_0 = \begin{pmatrix} u_0 & v_0 \\ u_0^{\prime} & v_0^{\prime} \end{pmatrix}.
\]
Then 
\[
\frac{dM}{dx} = \begin{pmatrix} 0 & 1 \\ V-z & 0 \end{pmatrix} M, \qquad
\frac{dM_0}{dx} = \begin{pmatrix} 0 & 1 \\ -z & 0 \end{pmatrix} M_0,
\]
and we have
\[
\frac{d}{dx} \left( M_0^{-1} M \right) = V A \left( M_0^{-1} M \right),
\]
where
\[
A = A(x, z) = \left( -v_0, u_0 \right)^T \left( u_0, v_0 \right)
\]
and $V=V(x)$. Hence
\[
\left( M_0^{-1} M \right) (x) = I + \int_0^x V(t) A(t) \left( M_0^{-1} M \right) (t) \, dt,
\]
where $I$ is the $2 \times 2$ identity matrix and, for notational convenience, we have suppressed 
the dependence on $z$. If $\Vert A \Vert$ denotes operator norm of the matrix $A$ in the two-dimensional 
space $l_2$, we have, for $x \ge 0$,
\[
\left\Vert \left( M_0^{-1} M \right) (x) - I \right\Vert \le 
\int_0^x \left| V(t) \right| \, \Vert A(t) \Vert \, dt + 
\int_0^x \left| V(t) \right| \, \Vert A(t) \Vert \, 
\left\Vert \left( M_0^{-1} M \right) (t) - I \right\Vert\, dt.
\]
An application of the Gronwall inequality now leads to the bound, valid for all $x \in \left[ 0, N \right]$,
\begin{align}
\Vert M(x) - M_0(x) \Vert 
&\le \Vert M_0(x) \Vert \, \left\Vert \left( M_0^{-1} M \right) (x) - I \right\Vert \nonumber\\
&\le \Vert M_0(x) \Vert 
\left\{ \exp \left( \int_0^N \left| V(t) \right| \Vert A(t) \Vert \, dt \right) - 1 \right\}.
\label{eminus1}
\end{align}
Noting that 
\[
\Vert M_0(x) \Vert \le \left( |u_0|^2 + |v_0|^2 + |u_0^{\prime}|^2 + |v_0^{\prime}|^2 \right)^{1/2},
\]
and 
\[
\Vert A \Vert = |u_0|^2 + |v_0|^2,
\]
we see that both $\Vert M_0(x, z)\Vert$ and $\Vert A(t, z)\Vert$ 
are bounded for $x, t \in \left[ 0, N \right]$ and $z \in K$. 

The result of the Lemma now follows from~\eqref{eminus1} and the observation that 
\[
\left| u-u_0 \right| \le \Vert M - M_0 \Vert, \qquad 
\left| v-v_0 \right| \le \Vert M - M_0 \Vert.
\]
\end{proof}

The following Corollary is a straightforward consequence of the Lemma.

\begin{corollary}
Let $K$ be a fixed compact subset of $\mathbb{C}^+$, and let $N>0$ be fixed. 
Define $u, v, u_0, v_0$ as in Lemma~\ref{lem:4}. 
Then given any $\epsilon > 0$, there exists $\delta_0 > 0$ 
($\delta_0$ depending on $\epsilon, N$ and $K$) 
such that, for all potential functions $V$ satisfying 
$\int_0^N \left| V(t) \right| \, dt < \delta_0$, 
we have, for all $z \in K$,
\begin{equation}
\left| 
\int_0^N \im \left( \vphantom{|^|} {\overline u}(t, z) v(t, z) \right) \, dt - 
\int_0^N \im \left( \vphantom{|^|} {\overline{u_0}}(t, z) v_0(t, z) \right) \, dt 
\right| 
< \epsilon.
\label{ubarv}
\end{equation}
\end{corollary}

\section{Estimate of $-f^{\prime}/f$ for potentials subject to an $L^2$-type condition}
We can now state an estimate of convergence of $-f^{\prime}/f$ to $i \sqrt{z}$ based 
on an $L^2$-type condition on the potential.

\begin{theorem}
\label{thm:1}
Let $f(x, z)$ be any solution for $x \in \left[ 0, \infty \right)$ 
of the Schr{\"o}dinger equation~\eqref{Schrodinger} at complex spectral parameter 
$z$ ($\im z > 0$) which satisfies the condition
\[
\im \left( - \frac{f^{\prime}(0, z)}{f(0, z)} \right) > 0.
\]
Let $K$ be any fixed compact subset of $\mathbb{C}^+$. 

Then, given any $\epsilon > 0$, there exist $\delta, N > 0$ 
($\delta, N$ depending on $\epsilon$ and $K$) such that, 
for all $L \ge N$ and for all potential functions $V$ satisfying 
the $L^2$ bound
\[
\int_0^L \left| V(t) \right|^2 \, dt < \delta,
\]
the estimate
\begin{equation}
\gamma \left( - \frac{f^{\prime}(L, z)}{f(L, z)}, i \sqrt{z} \right) < \epsilon
\label{festimate}
\end{equation}
holds for all $z \in K$.
\end{theorem}
\begin{proof}
In using the $\gamma$ measure of separation to carry out the estimate~\eqref{festimate}, 
it should be noted that, unlike the hyperbolic metric which is a function of $\gamma$, 
the separation $\gamma(z_1, z_2)$ between two points $z_1, z_2 \in \mathbb{C}^+$ 
does not satisfy the triangle inequality. However, the following result can be 
useful as a substitute for the triangle inequality: 

If $z_1, z_2, z_3 \in \mathbb{C}^+$ and it is given that
\[
\gamma(z_1, z_2) < \alpha, \qquad
\gamma(z_2, z_3) < \beta, \qquad 
\mbox{ with } \qquad 0 < \alpha, \beta \le 2,
\]
then it follows that $\gamma(z_1, z_3) < \sqrt{2} \left( \alpha + \beta \right)$. 
(To verify this result, note that if $0 < \alpha, \beta \le 2$ and 
\[
\gamma(z_1, z_2) = 2 \sinh \left( \frac{D(z_1, z_2)}{2} \right) < \alpha, \qquad
\gamma(z_2, z_3) = 2 \sinh \left( \frac{D(z_2, z_3)}{2} \right) <  \beta,
\]
then 
\begin{align*}
\gamma(z_1, z_3) &= 2 \sinh \left( \frac{D(z_1, z_3)}{2} \right)\\
&\le 2 \sinh \left( \frac{D(z_1, z_2)}{2} + \frac{D(z_2, z_3)}{2} \right)\\
&= 2 \sinh \left( \frac{D(z_1, z_2)}{2} \right) \cosh \left( \frac{D(z_2, z_3)}{2} \right) + 
2 \sinh \left( \frac{D(z_2, z_3)}{2} \right) \cosh \left( \frac{D(z_1, z_2)}{2} \right)\\
&\le \alpha \sqrt{1 + \frac{\beta^2}{4}} + \beta \sqrt{1 + \frac{\alpha^2}{4}}\\
&\le \left( \alpha + \beta \right) \sqrt{2}
\end{align*}
as required.) As a simple consequence of this result, the three inequalities 
$\gamma(z_1, z_2) < \frac{\epsilon}{6}, 
\gamma(z_2, z_3) < \frac{\epsilon}{6}, 
\gamma(z_3, z_4) < \frac{\epsilon}{6}$, 
with $0 < \epsilon < 1$, together imply that 
$\gamma(z_1, z_4) < \epsilon$. 

If, then, we define $u, v, u_0, v_0$ as in the proofs of the previous Lemmas,  
it will be sufficient, to verify~\eqref{festimate}, to show that if $z \in K$ 
then we have the three inequalities, at $x=L$,
\begin{equation}
\gamma \left( - \frac{f^{\prime}}{f}, - \frac{v^{\prime}}{v} \right) < \frac{\epsilon}{6}, \qquad
\gamma \left( - \frac{v^{\prime}}{v}, - \frac{v_0^{\prime}}{v_0} \right) < \frac{\epsilon}{6}, \qquad
\gamma \left( - \frac{v_0^{\prime}}{v_0}, i \sqrt{z} \right) < \frac{\epsilon}{6}.
\label{epsilon/6}
\end{equation}
We begin by fixing the value of $N$. Given $\epsilon > 0$ and a compact subset $K$ of $\mathbb{C}^+$, 
we take $N = N(\epsilon, K)$ to satisfy, for all $z \in K$, the three inequalities
\begin{subequations}
\begin{align}
\int_0^N \im \left( \vphantom{|^|} {\overline{u_0}} v_0 \right) \, dt &> \frac{12}{\epsilon \, \im z},
\label{inequality1}
\\
\frac{C^{\prime} \left( 1 + \left( \frac{b}{a} \right)^2 \right)^{1/2}}
{\left( \mathrm{e}^{4bN} - 1 \right)^{1/2}} &< \frac{\epsilon}{6},
\label{inequality2}
\\
N &> \frac{1}{\sqrt{|z|}}.
\label{inequality3}
\end{align}
\end{subequations}
That $N$ may be chosen to satisfy the first of these inequalities for $z \in K$ follows from the fact 
that $\int_0^{\infty} \im \left( {\overline{u_0}} v_0 \right) \, dt = \infty$ and that, 
for fixed $N$, the integral $\int_0^N \im \left( {\overline{u_0}} v_0 \right) \, dt$ 
depends continuously on $z$ for $\im z > 0$. In the second inequality we have 
$\sqrt{z} = a + i b$, where both $b$ and $b/a$ are bounded for $z \in K$; 
the constant $C^{\prime}$ is  defined in the proof of Lemma~\ref{lem:3}. 
Note also that $1/\sqrt{|z|}$ is bounded for $z \in K$ in the third inequality.

{}From the Corollary to Lemma~\ref{lem:4} we know that, for $z \in K$, 
the integral $\int_0^N \im \left( \vphantom{|^|} {\overline{u}} v \right) \! dt$ is close to 
$\int_0^N \im \left( \vphantom{|^|} {\overline{u_0}} v_0 \right) dt$ provided that 
$\int_0^N \left| V(t) \right| \, dt$ is sufficiently small. 
In particular, the inequality~\eqref{inequality1} implies that there exists 
$\delta_0 = \delta_0(\epsilon, K) > 0$ such that, for all $z \in K$, we have
\begin{equation}
\int_0^N \left| V(t) \right| \, dt < \delta_0 \Rightarrow 
\int_0^N \im \left( {\overline{u}} v \right) \, dt > \frac{6}{\epsilon \, \im z}.
\label{L21}
\end{equation}
Having fixed the values of $N$ and $\delta_0$, now define $\delta = \delta(\epsilon, K)$ 
to satisfy the two inequalities
\begin{itemize}
\item[(i)] $N \delta < \delta_0^2$;
\item[(ii)] $\frac{C |z|^{1/4}}{\im z} \sqrt{\delta} < \frac{\epsilon}{6}$ for all $z \in K$.
\end{itemize}
Here the constant $C$ has been defined in the statement of Lemma~\ref{lem:2}. 
Now suppose that $L \ge N$ and $\int_0^L \left| V(t) \right|^2 \, dt < \delta$. 
By the Schwarz inequality we then have
\[
\int_0^N \left| V(t) \right| \, dt \le 
\left( N \int_0^N \left| V(t) \right|^2 \, dt \right)^{1/2} < 
\left( \delta N \right)^{1/2} < \delta_0,
\]
by inequality~(i). Hence, \eqref{L21} implies that
\[
\int_0^N \im \left( {\overline{u}} v \right) \, dt > \frac{6}{\epsilon \, \im z}.
\]
By Lemma~3 of \cite{BP1} (see also Lemma~2 of \cite{BP2}) we have, for any solution 
$f$ of~\eqref{Schrodinger} satisfying $\im \left( - f^{\prime}(0, z) / f(0, z) \right) > 0$, 
\[
\left. \gamma \left( - \frac{f^{\prime}}{f}, - \frac{v^{\prime}}{v} \right) \right|_{x=L} \le 
\frac{1}{\im z  \int_0^L \im \left( {\overline{u}} v \right) \, dt} < \frac{\epsilon}{6}.
\]
Thus we have derived the first inequality in~\eqref{epsilon/6}. The second inequality 
in~\eqref{epsilon/6} follows from~\eqref{z1/4estimate} and (ii) above, 
using $\int_0^L \left| V(t) \right|^2 \, dt < \delta$. 
We can also use Lemma~\ref{lem:3} with the inequality~\eqref{inequality2} 
to complete the proof of~\eqref{epsilon/6}, which also completes the proof of the Theorem.
\end{proof}

We now explore some consequences of Theorem~\ref{thm:1} in the case of $L^2$-sparse potentials. 
Let $V$ be an $L^2$-sparse potential. Then a sequence of subintervals 
$\{ ( a_k, b_k ) \}$ ($k=1, 2, 3, \dots$) of $\mathbb{R}^+$ can be found such that, 
with $L_k = b_k - a_k$, 
\[
\lim_{k \rightarrow \infty} L_k = \infty 
\qquad \mbox{and} \qquad
\lim_{k \rightarrow \infty} \int_{a_k}^{b_k} \left( V(t) \right)^2 dt = 0.
\]
Given a fixed, bounded, measurable subset $A$ of $\mathbb{R}$, 
having closure ${\overline A}$, and given any $\epsilon > 0$, 
we first of all find $d>0$ ($d$ depending on $\epsilon$ and $A$) 
such that $E_A(d) < \epsilon \left| A \right| / 2$. 
Here $E_A(\cdot)$ is the error estimate on the right hand side of~\eqref{vddiffest}, 
and from~\eqref{vddiffest} we deduce that
\begin{equation}
\left| \mathcal{M}_1(A, S) - \mathcal{M}_2(A, S) \right| < 2 \epsilon \left| A \right|,
\label{vddiff2}
\end{equation}
provided $\gamma \left( F_1 (z), F_2(z) \right) < \epsilon$ for all $z \in K$, 
where $K$ is the compact subset of $\mathbb{C}^+$ defined by the  conditions 
$\im z = d, \re z \in {\overline A}$. 

Now use Theorem~\ref{thm:1} to define $\delta$ and $N$ such that, 
for all $L \ge N$ and for all potentials $V$ satisfying the bound 
$\int_0^L \left| V(t) \right|^2 dt < \delta$ we have
\begin{equation}
\gamma \left( - \frac{f^{\prime}(L, z)}{f(L, z)}, i \sqrt{z} \right) < \epsilon.
\label{delNdefn}
\end{equation}
Here $f(\cdot, z)$ is a solution of the Schr{\"o}dinger equation~\eqref{Schrodinger} for which 
\[
\im \left( - \frac{f^{\prime}(0, z)}{f(0, z)} \right) > 0.
\]
We take $k$ sufficiently large (say $k>k_0$) so that $L_k \ge N$ and such that 
the bound $\int_{a_k}^{b_k} \left| V(t) \right|^2 dt < \delta$ is satisfied by 
our sparse potential $V$. 

We can now apply~\eqref{delNdefn} with $L=L_k$, where $f$ is a suitably chosen 
solution of the Schr{\"o}dinger equation~\eqref{Schrodinger}, but with potential modified 
by an appropriate change of $x$-coordinate. There are two separate cases to be 
considered: 

Firstly, define $f(x, z) = v(x+a_k, z)$ (for $0 \le x \le L_k = b_k - a_k$). 
Then, for $x \in \left[ 0, L_k \right]$, $f(\cdot, z)$ satisfies the 
Schr{\"o}dinger equation~\eqref{Schrodinger} with potential $V(x+a_k)$. 
Moreover, we have
\[
\int_0^{L_k} \left( V(t+a_k) \right)^2 dt = \int_{a_k}^{b_k} \left( V(t) \right)^2 dt < \delta.
\]
Hence~\eqref{delNdefn} is satisfied in this case, and we have
\[
\gamma \left( - \frac{v^{\prime}(b_k, z)}{v(b_k, z)}, i \sqrt{z} \right) < \epsilon.
\]
{}From~\eqref{vddiff2} we now deduce that the respective value distributions for the 
Herglotz functions $- v^{\prime}(b_k, z) / v(b_k, z)$ and $i \sqrt{z}$ 
differ by at most $2 \epsilon \left| A \right|$, for all $k>k_0$.

Secondly, let $F(\cdot, z)$ be a (non-trivial) solution in $L^2 \left( 0, \infty \right)$ 
of the Schr{\"o}dinger equation~\eqref{Schrodinger}, with sparse potential $V$. The $m$-function 
$m^{a_k}(z)$ for the Schr{\"o}dinger operator $-\frac{d^2}{dx^2} + V$ acting in $L^2 \left( a_k, \infty \right)$ 
is then given by 
\[
m^{a_k}(z) = \frac{F^{\prime}(a_k, z)}{F(a_k, z)}.
\]
We can now define $f(\cdot, z)$ by
\[
f(x, z) = F(b_k - x, z) \qquad (0 \le x \le L_k)
\]
so that $f(\cdot, z)$ satisfies the Schr{\"o}dinger equation with potential $V(b_k - x)$. 
Since $F^{\prime}(b_k, z) / F(b_k, z)$ has positive imaginary part, 
we also have $\im \left( -f^{\prime}(0, z) / f(0, z) \right) > 0$. 
In this case, an application of~\eqref{festimate} with $L=L_k$ results in 
the estimate
\[
\gamma \left( m^{a_k}(z), i \sqrt{z} \right) < \epsilon, 
\]
and it follows as before that the respective value distributions for the 
Herglotz functions $m^{a_k}$ and $i\sqrt{z}$ differ by at most $2 \epsilon \left| A \right|$, 
for all $k>k_0$.

The following Theorem summarises the situation regarding asymptotic value distribution 
in the case of $L^2$-sparse potentials%
\footnote{We are indebted to A.\ Pushnitski for pointing out the close connection 
between estimates of $m$-functions at complex $z$ and asymptotic resolvent estimates 
in the case of potentials with an $L^2$ condition.}. %
The Theorem implies in particular, for the special case of $L^2$ potentials, 
that the value distribution of $v^{\prime}(N, \lambda) / v(N, \lambda)$ 
approaches an asymptotic limit as $N \rightarrow \infty$.

\begin{theorem}
\label{thm:2}
Let $v(\cdot, \lambda)$ be the solution of the Schr{\"o}dinger equation at real spectral parameter 
$\lambda$, subject to initial conditions $v(0, \lambda) = 0, v^{\prime}(0, \lambda) = 1$, 
in the case of an $L^2$-sparse potential $V$.

Let $\{ \left( a_k, b_k \right) \}$ be a sequence of subintervals of $\mathbb{R}^+$, 
for which $\lim_{k \rightarrow \infty} \left( b_k - a_k \right) = \infty$ and 
$\lim_{k \rightarrow \infty} \int_{a_k}^{b_k} \left| V(t) \right|^2 dt = 0$. 

Then for Borel subsets $A, S$ of $\mathbb{R}$, with $\left| A \right| < \infty$, we have
\begin{gather*}
\lim_{k \rightarrow \infty} \frac{1}{\pi} \int_A \theta \left( m_+^{a_k}(\lambda), S \right) d\lambda = 
\frac{1}{\pi} \int_A \theta \left( i \sqrt{\lambda}, S \right) d\lambda,
\\
\lim_{k \rightarrow \infty} 
\left| 
\left\{ 
\lambda \in A : \frac{v^{\prime}(b_k, \lambda)}{v(b_k, \lambda)} \in S 
\right\}
\right| = 
\frac{1}{\pi} \int_A \theta \left( i \sqrt{\lambda}, - S \right) d\lambda.
\end{gather*}
\end{theorem}

The conclusion of the Theorem, 
which applies in the first instance in the case that $A$ is bounded and of finite measure, 
may be extended to the more general case in which $A$ is not necessarily bounded. 
(Let $A$ have finite measure. Given $\epsilon > 0$, fix $N$ sufficiently large 
that the complement of $\left[ -N, N \right] \cap A$ has measure less than $\epsilon$. 
Denoting by $A_N$ this truncated set, the theorem may be applied first of all to $A_N$, 
which is bounded. Since the integrals to be estimated are then within $\epsilon$ of the 
corresponding integrals for the set $A$, the more general conclusion follows on letting 
$\epsilon$ approach zero.)

\section{Spectral analysis}
\label{sec:spectral}

Here we present some consequences of Theorem~\ref{thm:2} for the 
spectral theory of Schr{\"o}dinger operators with $L^2$-sparse potentials. 
The first result implies that absolutely continuous spectrum can 
occur only for $\lambda > 0$. 

\begin{corollary}
\label{cor:2}
Suppose $V$ is $L^2$-sparse. Then the support of the a.c. measure $\mu_{ac}$ of 
$T = - \frac{d^2}{dx^2} + V$ is contained in $\mathbb{R}^+$.
\end{corollary}
\begin{proof}
Suppose the contrary. Then if $\mu_{ac}$ is the a.c. part of the spectral measure, 
we can find a subset $A$ of $\mathbb{R}^-$ having finite Lebesgue measure for 
which $\mu_{ac}(A) > 0$. Then $\left| A \right| > 0$, and we may also suppose 
that $A$ is a subset of an essential support of $\mu_{ac}$. 

Now define intervals $\left( a_k, b_k \right)$ as in Theorem~\ref{thm:2}, and 
set $N_k = \left( a_k + b_k \right) / 2$. Then $N_k$ may be regarded 
either as the left hand endpoint of an interval $\left( N_k, b_k \right)$, 
or as the right hand endpoint of an interval $\left( a_k, N_k \right)$. 
An application of Theorem~\ref{thm:2} then implies that
\begin{equation}
\lim_{k \rightarrow \infty} \frac{1}{\pi} 
\int_A \theta \left( m_+^{N_k}(\lambda), S \right) d\lambda = 
\frac{1}{\pi} \int_A \theta \left( i \sqrt{\lambda}, S \right) d\lambda,
\label{specest1}
\end{equation}
whereas
\begin{equation}
\lim_{k \rightarrow \infty} 
\left| 
\left\{ 
\lambda \in A : \frac{v^{\prime}(N_k, \lambda)}{v(N_k, \lambda)} \in S 
\right\}
\right| = 
\frac{1}{\pi} \int_A \theta \left( i \sqrt{\lambda}, - S \right) d\lambda.
\label{specest2}
\end{equation}

Since $A$ is a subset of an essential support of $\mu_{ac}$, we also have
\begin{equation}
\lim_{k \rightarrow \infty} 
\left[ \left| \left\{ 
\lambda \in A : \frac{v^{\prime}(N_k, \lambda)}{v(N_k, \lambda)} \in S 
\right\} \right| 
- 
\frac{1}{\pi} 
\int_A \theta \left( m_+^{N_k}(\lambda), S \right) d\lambda 
\right] = 0.
\label{esssup}
\end{equation}
(For a proof of this result, which holds for any sequence $N_k$ with $N_k \rightarrow \infty$, 
and for arbitrary locally $L^1$ potentials, see \cite{BP1}.) 
Equations~\eqref{specest1}, \eqref{specest2} and~\eqref{esssup} now imply that
\begin{equation}
\int_A \theta \left( i \sqrt{\lambda}, S \right) d\lambda = 
\int_A \theta \left( i \sqrt{\lambda}, - S \right) d\lambda.
\label{plusminus}
\end{equation}
However $i \sqrt{\lambda} \in \mathbb{R}^-$ for $\lambda \in A$, and taking $S = \mathbb{R}^-$ 
we see that the left-hand-side of~\eqref{plusminus} is strictly positive, 
whereas the right-hand-side is zero. 

Hence we have a contradiction, and the Corollary is proved.
\end{proof}

There are interesting applications of Corollary~\ref{cor:2} to $L^2$ perturbations 
of slowly oscillating potentials such as $\cos \sqrt{x}$. For example, if 
$V(x) = \cos \sqrt{x} + V_0$ with $V_0 \in L^2 \left( \mathbb{R}^+ \right)$, 
then $V(x)-1$ is an $L^2$-sparse potential, and it follows from Corollary~\ref{cor:2} 
that $T=-\frac{d^2}{dx^2}+V$ has no a.c. measure for $\lambda < 1$. 
(In fact, $\left[ -1, 1 \right]$ is contained in the singular spectrum of $T$; 
for related results on spectral theory with slowly oscillating potentials see \cite{S}.) 

We can also consider various perturbations of $L^2$-sparse potentials. 
A typical result is the following:

\begin{corollary}
\label{cor:3}
Let $V$ be a $L^2$-sparse potential. 
Define intervals $\{ \left( a_k, b_k \right) \}$, 
with $N_k = \left( a_k + b_k \right) / 2$, 
as in the proof of Corollary~\ref{cor:2}. 
Then the Schr{\"o}dinger operator 
\[
- \frac{d^2}{dx^2} + V(x) + \sum_{k=1}^{\infty} \delta \left( x - N_k \right)
\]
has purely singular spectral measure.
\end{corollary}
\begin{proof}
The proof follows from Theorem~\ref{thm:2}, using similar arguments to those applied in 
\cite{BP1,BP2} to the special case in which $V$ is a sparse rather than $L^2$-sparse potential.
\end{proof}

\begin{ack}
This work has been partially supported by the EPSRC.
\end{ack}


\begin{thebibliography}{9999}

\bibitem[BP1]{BP1} Breimesser S.V., Pearson D.B.:
Asymptotic value distribution for solutions of the Schr{\"o}dinger equation, 
\textit{Mathematical Physics, Analysis and Geometry} \textbf{3}, 
385--403 (2000).

\bibitem[BP2]{BP2} Breimesser S.V., Pearson D.B.:
Geometrical aspects of spectral theory and value distribution for Herglotz functions.
Preprint 2001.

\bibitem[CL]{CL} Coddington E.A., Levinson N.: 
\textit{Theory of Ordinary Differential Equations} 
(McGraw-Hill, New York, 1955). 

\bibitem[DK]{DK} Deift R., Killip R.:
On the absolutely continuous spectrum of one-dimensional
Schroedinger operators with square summable potentials,
\textit{Communications in Mathematical Physics} \textbf{203},
341--347 (1999).

\bibitem[G]{G} Glazman I.M.: 
\textit{Direct Methods of Qualitative Spectral Analysis of Singular Differential Operators} 
(Israel program for scientific translations, Jerusalem, 1965).

\bibitem[KLS]{KLS} Kiselev A., Last Y., Simon B.: 
Modified Pr{\"u}fer and EFGP transforms and the spectral analysis of one-dimensional Schr{\"o}dinger operators, 
\textit{Communications in Mathematical Physics} \textbf{194}, 
1--45 (1998).

\bibitem[R]{R} Remling C.: 
A probabilistic approach to one-dimensional Schr{\"o}dinger operators with sparse potentials, 
\textit{Communications in Mathematical Physics} \textbf{185}, 
313--323 (1997).

\bibitem[SS]{SS} Simon B., Stolz G.: 
Operators with singular continuous spectrum, V. Sparse potentials, 
\textit{Proceedings of the American Mathematical Society} \textbf{124}, 
2073--2080 (1996).

\bibitem[S]{S} Stolz G.: 
Spectral theory for slowly oscillating potentials II. 
Schr{\"o}dinger operators, 
\textit{Mathematische Nachrichten} \textbf{183}, 
275--294 (1997).

\end{thebibliography}
\end{document}